\documentclass{article}
\usepackage{amsthm,amsmath,amssymb,hyperref}

\newcommand{\nc}{\newcommand}

\newtheorem{thm}{Theorem}[section]
\newtheorem{rmk}[thm]{Remark}

\newtheorem{prop}[thm]{Proposition}
\newtheorem{lemma}[thm]{Lemma}

\newtheorem{corollary}[thm]{Corollary}

\newtheorem{definition}[thm]{Definition}

\newenvironment{defin}{\begin{definition} \rm}{\end{definition}}
\newenvironment{cor}{\begin{corollary} \rm}{\end{corollary}}
\newenvironment{lem}{\begin{lemma}\rm }{\end{lemma}}

\nc{\Ext}{\operatorname{Ext}}
\nc{\NS}{\operatorname{NS}}
\nc{\Amp}{\operatorname{Amp}}
\nc{\Pic}{\operatorname{Pic}}
\nc{\Kom}{\operatorname{Kom}}
\nc{\gr}{\operatorname{gr}}
\nc{\Gr}{\operatorname{Gr}}
\nc{\Rep}{\operatorname{Rep}}
\nc{\Hom}{\operatorname{Hom}}
\nc{\RHom}{R\operatorname{Hom}}
\nc{\cRHom}{\operatorname{\mathcal{R}\mathcal{H}om}}
\nc{\cHom}{\operatorname{\mathcal{H}om}}
\nc{\End}{\operatorname{End}}
\nc{\Coh}{\operatorname{Coh}}
\nc{\Aut}{\operatorname{Aut}}
\nc{\Coker}{\operatorname{Coker}}
\nc{\coker}{\operatorname{coker}}
\nc{\Ker}{\operatorname{Ker}}
\nc{\Obj}{\operatorname{Obj}}
\nc{\img}{\operatorname{Im}}
\nc{\D}{\operatorname{D}}
\nc{\ch}{\operatorname{ch}}
\nc{\Stab}{\operatorname{Stab}}
\nc{\SL}{\operatorname{Stab}}
\nc{\rk}{\operatorname{rk}}
\nc{\GL}{\operatorname{GL}}
\nc{\Log}{\mathop{\mathrm{Log}}}
\nc{\abs}[1]{\lvert#1\rvert}
\nc{\Cone}{\operatorname{Cone}}
\nc{\id}{\operatorname{id}}

\nc{\cA}{{\mathcal A}}
\nc{\cB}{{\mathcal B}}
\nc{\cC}{{\mathcal C}}
\nc{\cD}{{\mathcal D}}
\nc{\cE}{{\mathcal E}}
\nc{\cF}{{\mathcal F}}
\nc{\cG}{{\mathcal G}}
\nc{\cH}{{\mathcal H}}
\nc{\cI}{{\mathcal I}}
\nc{\cJ}{{\mathcal J}}
\nc{\cK}{{\mathcal K}}
\nc{\cL}{{\mathcal L}}
\nc{\cM}{{\mathcal M}}
\nc{\cN}{{\mathcal N}}
\nc{\cO}{{\mathcal O}}
\nc{\cP}{{\mathcal P}}
\nc{\cQ}{{\mathcal Q}}
\nc{\cR}{{\mathcal R}}
\nc{\cS}{{\mathcal S}}
\nc{\cT}{{\mathcal T}}
\nc{\cU}{{\mathcal U}}
\nc{\cV}{{\mathcal V}}
\nc{\cW}{{\mathcal W}}
\nc{\cX}{{\mathcal X}}
\nc{\cY}{{\mathcal Y}}
\nc{\cZ}{{\mathcal Z}}

\nc{\bA}{{\mathbb A}}
\nc{\bB}{{\mathbb B}}
\nc{\bC}{{\mathbb C}}
\nc{\bD}{{\mathbb D}}
\nc{\bE}{{\mathbb E}}
\nc{\bF}{{\mathbb F}}
\nc{\bG}{{\mathbb G}}
\nc{\bH}{{\mathbb H}}
\nc{\bI}{{\mathbb I}}
\nc{\bJ}{{\mathbb J}}
\nc{\bK}{{\mathbb K}}
\nc{\bL}{{\mathbb L}}
\nc{\bM}{{\mathbb M}}
\nc{\bN}{{\mathbb N}}
\nc{\bO}{{\mathbb O}}
\nc{\bP}{{\mathbb P}}
\nc{\bQ}{{\mathbb Q}}
\nc{\bR}{{\mathbb R}}
\nc{\bS}{{\mathbb S}}
\nc{\bT}{{\mathbb T}}
\nc{\bU}{{\mathbb U}}
\nc{\bV}{{\mathbb V}}
\nc{\bW}{{\mathbb W}}
\nc{\bX}{{\mathbb X}}
\nc{\bY}{{\mathbb Y}}
\nc{\bZ}{{\mathbb Z}}

\begin{document}
\title{Topologies on a triangulated category}
\author{So Okada\footnote{Address: Max-Planck-Institut f\"ur Mathematik,
Vivatsgasse 7, Bonn Germany 53111, Email: okada@mpim-bonn.mpg.de}}
\maketitle
\begin{abstract}
 On objects of a triangulated category with a {\it stability condition},
 we construct a topology.
\end{abstract}
\section{Introduction}

For any triangulated category, Bridgeland \cite{math.AG/0212237}
introduced the notions of {\it stability conditions} and {\it stability
manifolds}, motivated by the Douglas' work \cite{MR1840318},
\cite{MR1909947}, \cite{MR1957548} on $\Pi$-stabilities of D-branes for
Calabi-Yau manifolds.

 In particular, for a K3 surface $X$ and the bounded derived category of
 coherent sheaves on $X$, denoted by $D(X)$, Bridgeland
 \cite{math.AG/0307164} proved that stability conditions on the
 stability manifold approximate Gieseker stabilities \cite{MR466475},
 \cite{Mar77}, \cite{MR509499}.
 
 Roughly speaking in terms of string theory, for a Calabi-Yau manifold
 $X$, objects of $D(X)$ correspond to B-branes among D-branes, that are
 boundary conditions of open strings.  For a stability condition, among
 B-branes, {\it semistable} objects correspond to BPS states, that
 recognize themselves as D-branes in the untwisted topological field
 theory. A physical quantity of each BPS state is called a {\it central
 charge}. The author recommends \cite{Aspinwall:2004jr} for string
 theory related to this subject.

 For a triangulated category $\cT$, stability conditions on the
 stability manifold describe variation of $t$-structures of $\cT$, in
 collaboration with other subjects.  The author recommends
 \cite{math.AG/0611510}, \cite[Section 13]{MR2244106},
 \cite{math.AG/0602129} for introductions to this subject.

 In some general settings, for a triangulated category, we have seen
 spectra \cite[Section 12]{MPI-2003-112}, \cite{MR2196732} and moduli
 spaces of objects with vanishing conditions on extension groups
 \cite{MR1966840}, \cite{MR2177199}, \cite{math.AG/0612078}.

  In this article, on objects of a triangulated category with a
  stability condition, we construct a topology induced from the central
  charge.  Also, with a {\it faithful} or a {\it numerically faithful
  stability condition}, we find our topology compatible with the {\it
  Grothendieck group} or the {\it numerical Grothendieck group}.  We
  realize that objects of a triangulated category with a stability
  condition is a connected space.

\subsubsection*{Acknowledgments}
The author thanks I. Mirkovi\'c and T. Mochizuki for their
discussions. The author is grateful to the Max-Planck-Institut f\"{u}r
Mathematik for their excellent support and working environment and to
{\it Bonner} for their warm hospitality.

 \section{Stability conditions and stability manifolds}
 Throughout this paper, $\cT$ is a triangulated category such that for
 any objects $E$ and $F$ of $\cT$, $\oplus_{i\in
 \bZ}\Hom^{i}_{\cT}(E,F)$ is a finite-dimensional vector space over
 $\bC$.  For example, $\cT$ can be $D(X)$ of a smooth projective variety
 $X$ over $\bC$.

 In the {\it Grothendieck group} $K(\cT)$ of $\cT$ and for each object
 $E$ in $\cT$, let $[E]$ be the class of $E$.  For objects $E$ and $F$
 of $\cT$, the {\it Euler paring} $\chi(E,F)$ is defined to be
 $\sum_{i\in \bZ}(-1)^{i}\dim \Hom^{i}_{\cT}(E,F)$.  On this paring, the
 {\it numerical Grothendieck group} $N(\cT)$ is defined to be the
 quotient of $K(\cT)$ by $K(\cT)^{\perp}$.

 From \cite{math.AG/0212237}, we will recall fundamental notions.

\subsection{Stability conditions}
A {\it stability condition} $\sigma=(Z, \cP)$ on $\cT$ consists of a
group homomorphism from $K(\cT)$ to $\bC$, called a {\it central charge}
$Z$, and a family of full abelian subcategories of $\cT$, called a {\it
slicing} $\cP(k)$, indexed by real numbers $k$, with the following
conditions.

 For each real number $k$, if $E$ is an object of $\cP(k)$, then for
 some positive real number $m(E)$, called the {\it mass} of $E$, we have
 $Z(E)=m(E)\exp(i\pi k )$. For each real number $k$, we have
 $\cP(k+1)=\cP(k)[1]$. For any real numbers $k_{1}>k_{2}$ and any
 objects $E_{i}$ of $\cP(k_{i})$, \ $\Hom_{\cT}(E_{1}, E_{2})$ is the
 zero vector space. For any nonzero object $E$ of $\cT$, there exists a
 finite sequence of real numbers $k_{1}>\cdots>k_{n}$ and objects
 $H^{k_{i}}_{\sigma}(E)$ of $\cP(k_{i})$ such that there exists a
 sequence of exact triangles $E_{i-1}\to E_{i}\to H_{\sigma}^{k_{i}}(E)$
 with $E_{0}$ and $E$ being the zero object and $E$.

 The above sequence of the exact triangles is unique up to isomorphisms
 and called the {\it Harder-Narasimhan filtration} of $E$; also, we call
 a real number in the above sequence of the real numbers a {\it
 nontrivial phase} of $E$.  For each real number $k$, any nonzero object
 of $\cP(k)$ is called {\it semistable}.

 If the central charge $Z$ factors through $N(\cT)$, then $\sigma$ is
 called a {\it numerical stability condition}.

\subsection{Hearts of stability conditions} 
For an interval $I$ in real numbers, $\cP(I)$ is defined to be the
smallest full subcategory of $\cT$ consisting of objects of $\cP(k)$ for
each real number $k$ in $I$, it is closed under extension; i.e., if
$E\to G\to F$ is an exact triangle in $\cT$ and both $E$ and $F$ are
objects of $\cP(I)$, then $G$ is an object of $\cP(I)$.  In particular,
for each real number $j$, $\cP((j-1,j])$ is a heart of a bounded
$t$-structure of $\cT$.  We will call all $\cP((j-1,j])$ for real
numbers $j$, ``hearts of $\sigma$''.

  For each nonzero object $E$ of $\cP((j-1, j])$, the {\it phase} of $E$
  is defined to be $\phi(E)=(1/\pi)\arg Z(E)\in(j-1, j]$.

  \subsection{Stability manifolds}
  A subset of stability conditions on $\cT$ makes the {\it stability
  manifold} $\Stab(\cT)$, this has a natural topology induced from the
  central charges and each connected component is a manifold locally
  modeled on some topological vector subspace of $\Hom_{\bZ}(K(\cT),
  \bC)$.

   The subset of $\Stab(\cT)$ consisting of numerical stability
   conditions makes a subspace, the {\it numerical stability manifold}
   $\Stab_{N}(\cT)$, this is locally modeled on some topological vector
   subspace of $\Hom_{\bZ}(N(\cT), \bC)$.

 \section{Faithful or numerically faithful stability conditions}
 Let $K(\cT)_{\bQ}$ and $N(\cT)_{\bQ}$ denote the tensor products
 $K(\cT)\otimes \bQ$ and $N(\cT)\otimes \bQ$.

    \begin{defin}
     Let $\sigma$ be a stability condition on $\cT$.  We call $\sigma$
     {\it faithful}, if whenever nonzero objects $E$ and $F$ of a heart
     of $\sigma$ are linearly independent in $K(\cT)_{\bQ}$, then
     objects $E$ and $F$ are with different phases.  Likewise, we call
     $\sigma$ {\it numerically faithful}, if whenever nonzero objects
     $E$ and $F$ of a heart of $\sigma$ are linearly independent in
     $N(\cT)_{\bQ}$, then $E$ and $F$ are with different phases.
    \end{defin}
    
    In any stability manifolds that we are aware of, by the following
    lemma, there are faithful or numerically faithful stability
    conditions.
 
   \begin{lem}\label{lem:max}
    If $K(\cT)_{\bQ}$ has no more than countable dimension over $\bQ$,
    and if $\Stab(\cT)$ carries a connected component $M$ that is
    locally isomorphic to $\Hom_{\bZ}(K(\cT), \bC)$, then the subset of
    $M$ consisting of faithful stability conditions is dense in $M$.
    Likewise, if $N(\cT)_{\bQ}$ has no more than countable dimension
    over $\bQ$, and if $\Stab_{N}(\cT)$ carries a connected component
    $M$ that is locally isomorphic to $\Hom_{\bZ}(N(\cT), \bC)$, then
    the subset of $M$ consisting of numerically faithful stability
    conditions is dense in $M$.
   \end{lem}
   \begin{proof}
    Let $T$ be the subset of $\Hom_{\bZ}(K(\cT), \bC)$ consisting of $Z$
    such that for some linearly independent classes $[E]$ and $[F]$ in
    $K(\cT)_{\bQ}$, in the interval $(0,2]$, real numbers $(1/\pi)\arg
    Z(E)$ and $(1/\pi)\arg Z(F)$ are the same. Then, $T$ is a countable
    union of codimension-one subspaces of $\Hom_{\bZ}(K(\cT), \bC)$.
    Thus, the complement of $T$ in $\Hom_{\bZ}(K(\cT), \bC)$ is dense in
    $\Hom_{\bZ}(K(\cT), \bC)$.  The same argument holds for the latter
    case.
   \end{proof}
   
   In particular, we have the following.

   \begin{cor}\label{cor:Q}
    For a stability condition $\sigma$ of $\cT$, let $E$ and $F$ be
    nonzero objects of a heart of $\sigma$ with the same phases.  If
    $\sigma$ is a faithful stability condition, then for some positive
    rational number $q$, \ $[E]=q [F]$ in $K(\cT)_{\bQ}$. If $\sigma$ is
    a numerically faithful stability condition, then for some positive
    rational number $q$, \ $[E]=q [F]$ in $N(\cT)_{\bQ}$.
   \end{cor}
   \begin{proof}
    Since $\sigma=(Z, \cP)$ is faithful, $[E]$ and $[F]$ are linearly
    dependent in $K(\cT)_{\bQ}$. So there exists a nonzero rational
    number $q$ such that $[E]=q[F]$ in $K(\cT)_{\bQ}$.  Here, $q$ can
    not be zero, since $E$ and
    $F$ are nonzero objects of a heart of $\sigma$,
    which implies $Z(E)$ and $Z(F)$ are not zero.  Also, $q$ can not be
    negative; otherwise, their phases would differ by an odd
    integer. The same argument holds for the latter case.
   \end{proof}

 \section{Topologies on a triangulated category}
   We notice that for a stability condition $(Z, \cP)$ on $\cT$ and each
   semistable object of $\cT$, the central charge $Z$ factors through
   $\bP^{1}$; i.e., for each semistable object $E$ of $\cT$, its central
   charge is the exponential of the complex number
   $\log(m(E))+i\pi\phi(E)$ in $\bP^{1}$.  Now, we define the following.

 \begin{defin}\label{def:extended}
  Let $\sigma=(Z, \cP)$ be a stability condition on $\cT$.  For a
  nonzero object $E$ of $\cT$, let $\tilde{Z}(E)$ be the subset of
  $\bP^{1}$ consisting of points $\log(m(H_{\sigma}^{k}(E)))+i\pi
  \phi(H_{\sigma}^{k}(E))$ for nontrivial phases $k$ of $E$. For the
  zero object, let the image of $\tilde{Z}$ be the infinite point in
  $\bP^{1}$.  We call the function $\tilde{Z}$ the {\it extended central
  charge} of $\sigma$.
 \end{defin} 

 For our arguments here, let us assume the Euclidean topology on
 $\bP^{1}$.

 \begin{prop}\label{prop:topo}
  For a triangulated category with a stability condition, there exists a
  unique topology on objects of the triangulated category such that the
  extended central charge is continuous.
 \end{prop}
 \begin{proof}
  Let $\cT$ be a triangulated category with a stability condition
  $\sigma=(Z, \cP)$, and $\tilde{Z}$ be the extended central
  charge. Then, finite unions of $\tilde{Z}^{-1}(V)$ for closed subsets
  $V$ of $\bP^{1}$ are our closed subsets of objects of $\cT$.
 \end{proof}

 On the topology in Proposition \ref{prop:topo}, we have the following
 corollaries.

 \begin{cor}\label{cor:point_set}
  For a semistable object $E$ of $\cT$, any object of the closed set
  $\tilde{Z}^{-1}(\tilde{Z}(E))$ is semistable.
 \end{cor}
 \begin{proof}
  The Harder-Narasimhan filtration of any object of
  $\tilde{Z}^{-1}(\tilde{Z}(E))$ is trivial, since its image under
  $\tilde{Z}$ is a point in $\bP^{1}$.
 \end{proof}

 In particular, we have the following.

 \begin{cor}\label{cor:closed_set}
  For a faithful stability condition $(Z, \cP)$ on $\cT$ and a
  semistable object $E$ of $\cT$, any object of the closed set
  $\tilde{Z}^{-1}(\tilde{Z}(E))$ is semistable and has the class $[E]$
  in $K(\cT)$.  For a numerically faithful stability condition $(Z,
  \cP)$ on $\cT$ and a semistable object $E$ of $\cT$, any object of the
  closed set $\tilde{Z}^{-1}(\tilde{Z}(E))$ is semistable and has the
  class $[E]$ in $N(\cT)$.
\end{cor}
 \begin{proof}
  By Corollaries \ref{cor:Q} and \ref{cor:point_set}, these statements
  hold.
 \end{proof}

  \subsection{Connectedness}

 \begin{lem}\label{lem:zero}
  Let $\cT$ be a triangulated category with a stability condition.  On
  the topology in Proposition \ref{prop:topo}, if an open subset $U$ of
  objects of $\cT$ contains the zero object of $\cT$, then $U$ contains
  some semistable objects.
 \end{lem}
 \begin{proof}
  Let $(Z, \cP)$ denote the stability condition.  Let $\tilde{Z}$ be the
  extended central charge.  Now, the open subset $U$ is a complement of
  the union $\tilde{Z}^{-1}(V_{1})\cup \cdots \cup \tilde{Z}^{-1}(V_{n})$
  for some closed subsets $V_{i}$ of $\bP^{1}$.  Since the zero object is
  not in $U$, there is an open subset $U'$ of $\bP^{1}$ such that
  $\tilde{Z}^{-1}(U')$ is a subset of $U$ consisting of the infinite
  point of $\bP^{1}$.

  Here, since for any semistable object $E$ of $\cT$, the direct sum
  $E\oplus E$ is again semistable with the doubled mass, so
  $\tilde{Z}^{-1}(U')$ contains some semistable objects.
 \end{proof}

 \begin{lem}\label{lem:all}
  Let $\cT$ be a triangulated category with a stability condition.  On
  the topology constructed in Proposition \ref{prop:topo}, if a proper
  open subset $U$ of objects of $\cT$ contains all semistable objects of
  $\cT$, then the open set $U$ does not contain the zero object.
 \end{lem}
 \begin{proof}
  For closed subsets $V_{i}$ of $\bP^{1}$ such that
  $U=\tilde{Z}^{-1}(V_{1})\cup \cdots \cup \tilde{Z}^{-1}(V_{n})$, by
  the assumption on $U$, there is no semistable object $E$ of $\cT$ such
  that for some $V_{i}$, \ $\tilde{Z}(E)$ is in $V_{i}$.  So, by Definition
  \ref{def:extended}, either $\tilde{Z}^{-1}(V_{i})$ is empty or of only
  the zero object.
 \end{proof}
 
\begin{thm}\label{thm:conn}
 For any triangulated category with a stability condition, on the
 topology in Proposition \ref{prop:topo}, objects of the triangulated
 category is a connected space
\end{thm}
\begin{proof}
 Let $\cT$ be a triangulated category with a stability condition $(Z,
 \cP)$, and $\tilde{Z}$ be the extended central charge.  Let us prove
 the statement by a contradiction. So we suppose that some open subsets
 $U_{1}$ and $U_{2}$ of objects of $\cT$ separate objects of $\cT$.
 Then, either there exist semistable objects $E_{1}$ and $E_{2}$ of
 $\cT$ such that $E_{1}$ and $E_{2}$ are objects of $U_{1}$ and $U_{2}$,
 or one of $U_{1}$ and $U_{2}$ contains all semistable objects.  For the
 former case, the direct sum $E_{1}\oplus E_{2}$ is not a object of
 $U_{1}\cup U_{2}$.  For the latter case, if the open set $U_{1}$
 contains all semistable objects, then by Lemma \ref{lem:all}, $U_{2}$
 is an open subset consisting of the zero object. However, then by Lemma
 \ref{lem:zero}, $U_{1}$ and $U_{2}$ can not be disjoint.
\end{proof}

\begin{rmk}
 For any triangulated category with a stability condition, the extended
 central charge naturally factors through the {\it Ran's space}
 \cite[3.4]{MR2058353} of $\bP^{1}$ with respect to a topology on
 $\bP^{1}$.
\end{rmk}

  \subsection{Closed sets}
  For a faithful stability condition $\sigma=(Z, \cP)$ on $\cT$
  and some semistable objects $E$ of $\cT$, let us take a look at closed sets
  $\tilde{Z}^{-1}(\tilde{Z}(E))$ in Corollary \ref{cor:closed_set}.  For a
  vector space $V$, let $V^{*}$ denote the dual of $V$.  From here, let
  the shift $[2]$ be the {\it Serre functor} \cite[Definition 3.1]{BoKa}
  of $\cT$; i.e., for any objects $E$ and $F$ of $\cT$, there exist
  bifunctorial isomorphisms $\Psi_{E,F}:\Hom_{\cT}(E,F)\cong
  \Hom_{\cT}(F, E[2])^{*}$ such that $(\Psi^{-1}_{E[2],F[2]})^{*}\circ
  \Psi_{E,F}$ coincides with the isomorphism induced by $[2]$.  In other
  words, $\cT$ is a {\it 2-Calabi-Yau category} \cite{Kon}.

  If an object $E$ of $\cT$ satisfies that for any integer $i$ other
  than $0$ or $2$, \ $\Hom_{\cT}^{i}(E,E)$ is the zero vector space and
  $\Hom_{\cT}(E,E)\cong \Hom_{\cT}^{2}(E,E)^{*}=\bC$, then $E$ is called
  {\it spherical} \cite[Definition 1.1]{SeiTho}.

  For semistable objects $E$ and $E'$ of $\cT$ with the same phases,
  since they are of a heart of $\cT$, for 
  any negative integer $i$, 
  $\Hom_{\cT}^{i}(E,E')$ is the zero vector space, and since
  $\Hom_{\cT}^{i}(E,E')\cong \Hom_{\cT}^{2-i}(E',E)^{*}$, for 
  any integer $i>2$,  \
  $\Hom_{\cT}^{i}(E,E')$ is the zero vector space.

  Let us recall that for a stability condition $(Z, \cP)$ on $\cT$, a
  semistable object $E$ of $\cP(\phi(E))$ is called {\it stable} if it
  has no nontrivial subobject of the abelian category $\cP(\phi(E))$.
  In particular, any stable object $E$ of $\cT$ satisfies
  $\Hom_{\cT}(E,E)=\bC$.  Since now $\cT$ is 2-Calabi-Yau, a simple case
  of the paring in \cite[Proposition I.1.4]{ReiVan} tells us that
  $\chi(E,E)$ is even \cite[Remark 3.7]{Oka}.  

  Each stability condition in $\Stab(\cT)$ has the property called {\it
  locally-finiteness}; in particular, for any locally-finite stability
  condition and a semistable object $E$, in $\cP(\phi(E))$, we have a
  Jordan-H\"older decomposition whose composition factors are stable
  objects, called {\it stable factors} of $E$.
  
  Then, in our formality, we realize a part of \cite[Corollary 3.6]{Mu}.

  \begin{prop}\label{prop:point}
   Let $\sigma=(Z, \cP)$ be a faithful or numerically faithful stability
   condition in $\Stab(\cT)$ and $E$ be a semistable object of $\cT$
   such that $\chi(E,E)$ is positive.  Then, $\tilde{Z}^{-1}(\tilde{Z}(E))$
   contains only $E$, and $E$ is a direct sum of a stable spherical
   object.
  \end{prop}
  \begin{proof}
   Any stable object $S$ of $\cP(\phi(E))$ is spherical; because,
   $\chi(S,S)=2- \Hom_{\cT}^{1}(S,S)$ is even and by Corollary
   \ref{cor:Q}, for some positive rational number $q$, $[S]=q[E]$, which
   implies $\chi(S,S)=q^{2}\chi(E,E)$ is positive.  If we have
   nonisomorphic stable objects $S_{1}$ and $S_{2}$ in $\cP(\phi(E))$,
   then $\Hom_{\cT}(S_{1},S_{2})$ and $\Hom_{\cT}^{2}(S_{1},S_{2}) \cong
   \Hom_{\cT}(S_{2},S_{1})^{*}$ would be the zero vector spaces, and then,
   $\chi(S_{1},S_{2})=- \dim\Hom^{1}_{\cT}(S_{1},S_{2})$ 
   would not be positive. However,
   since by Corollary \ref{cor:Q}, for some positive rational numbers
   $q_{1}$ and $q_{2}$, we have $[E]=q_{1}[S_{1}]=q_{2}[S_{2}]$,  
   \ $\chi(S_{1},S_{2})=q_{1}q_{2}\chi(E,E)$ is positive.

   Extensions of a spherical object are direct sums of the spherical
   object, so the statement follows.
  \end{proof}
  
  \subsubsection{An example}
  Let $X$ be the cotangent bundle of $\bP^{1}$. From here, let $\cT$ be
  the full subcategory of $D(X)$ consisting of objects supported over
  $\bP^{1}$.
  
  Since $K(\cT)\cong \bZ [\cO_{x}] \oplus \bZ [\cO_{\bP^{1}}]$ with
  $\chi(\cO_{\bP^{1}},\cO_{\bP^{1}})=2$ and $[\cO_{x}]$ being zero in
  $N(\cT)$, for any object $E$,  $\chi(E,E)$ is
  not negative.  So, by Proposition
  \ref{prop:point}, we take semistable objects $E$ with $[E]$ being zero
  in $N(\cT)$ and look at $\tilde{Z}^{-1}(\tilde{Z}(E))$.

  For a spherical object $E$ and an object $F$ of $\cT$, the cone of the
  evaluation map $\RHom_{\cT}(E, F)\otimes E \to F$ is denoted by
  $T_{E}(F)$, the {\it twist functor} of $E$ \cite[Section 1.1]{SeiTho}.
  By \cite[Theorem 1.2]{SeiTho}, twist functors are autoequivalences of
  $\cT$.  Any autoequivalence $\Phi$ of $\cT$ acts on $\Stab(\cT)$; for
  each $\sigma=(Z, \cP)$ in $\Stab(\cT)$, $\Psi(\sigma)$
  is defined to be $(Z\circ
  \Psi^{-1}, \Psi \circ \cP)$.

    Here, for an object $E$ of $\cT$, if for any point $x$ in $\bP^{1}$
    and any integer $i$, $\Hom_{\cT}^{i}(E, \cO_{x})$ is the zero vector
    space, then $E$ is isomorphic to the zero object; for the largest
    integer $q_{0}$ such that the support of the cohomology sheaf
    $H_{X}^{q_{0}}(E)$ has a point $x$ in $\bP^{1}$, the nonzero term
    $E_{2}^{0,q_{0}}$ in the spectral sequence $E_{2}^{p,q}
    =\Hom^{p}_{\cT}(H_{X}^{-q}(E), \cO_{x})\Rightarrow
    \Hom_{\cT}^{p+q}(E, \cO_{x})$ survives at infinity (see
    \cite[Section 2]{Bri99}).

  Let us recall that the connectedness of $\Stab(\cT)$ in \cite[Theorem
  4.12]{Oka} follows by proving that for any stability condition
  $\sigma$ in $\Stab(\cT)$ and for some integer $w$, objects
  $\cO_{\bP^{1}}(w-1)[1]$ and $\cO_{\bP^{1}}(w)$ of $\cT$ generate a
  heart of $\sigma$ and by using \cite[Lemmas 3.1 and 3.6, Theorem
  1.3]{Brid}, as explained in \cite[Section 4.3]{Oka}.

  Now, to prove Proposition \ref{prop:P1}, one way is to use
  \cite[Proposition 18 and Corollary 20]{IshUedUeh}; for our conclusion
  on this case, we give a proof.

  \begin{prop}\label{prop:P1}
   Let $\sigma=(Z, \cP)$ be a faithful stability condition in
   $\Stab(\cT)$.  For each semistable object $E$ of $\cT$ with
   $[E]=[\cO_{x}]$, up to autoequivalences of $\cT$,
   $\tilde{Z}^{-1}(\tilde{Z}(E))=\{ \cO_{x} \mid x \in \bP^{1}\}$.
  \end{prop}
  \begin{proof}
   As mentioned above, for some integer $w$, by extensions, objects
   $\cO_{\bP^{1}}(w-1)[1]$ and $\cO_{\bP^{1}}(w)$ generate a heart of
   $\sigma$.  For each point $x$ in $\bP^{1}$, by the exact triangle $
   \cO_{\bP^{1}}(w)\to \cO_{x} \to\cO_{\bP^{1}}(w-1)[1]$, 
   the object $\cO_{x}$ is
    of the heart and, since $\sigma$ is faithful, phases of
   $\cO_{\bP^{1}}(w-1)[1]$ and $\cO_{\bP^{1}}(w)$ are distinct.
   
   First, let us suppose $\phi(\cO_{\bP^{1}}(w-1)[1])>
   \phi(\cO_{\bP^{1}})$ and see for any point $x$ in $\bP^{1}$,
   $\cO_{x}$ is semistable. With classes of objects of the heart,
   $[\cO_{x}]$ can be represented by only
   $[\cO_{\bP^{1}}(w-1)[1]]+[\cO_{\bP^{1}}(w)]$, if it is other than
   $[\cO_{x}]$.  So, if $\cO_{x}$ were not semistable, by the assumption
   on phases, the Harder-Narasimhan filtration of $\cO_{x}$ must give
   the exact triangle $\cO_{\bP^{1}}(w-1)[1]\to \cO_{x} \to
   \cO_{\bP^{1}}(w)$; however, $\Hom_{\cT}(\cO_{\bP^{1}}(w-1)[1]),
   \cO_{\bP^{1}}(w))=\Hom_{\cT}^{-1}(\cO_{\bP^{1}}(w-1), \cO_{x})$ is
   the zero vector space, since the objects $\cO_{\bP^{1}}(w-1)$ and
   $\cO_{x}$ are of a heart of $\cT$ such as the category of the
   coherent sheaves of $X$ supported over $\bP^{1}$.

   Now, any object of $\tilde{Z}^{-1}(\tilde{Z}(E))$ is stable;
   otherwise, by Corollary \ref{cor:Q}, for some rational number $q>1$
   and a stable object $S$, we would have $[\cO_{x}]=q [S]$.

   For a point $x$ in $\bP^{1}$ and the integer
   $k=\phi(E)-\phi(\cO_{x})$, we show that for each object $E'$ of
   $\tilde{Z}^{-1}(\tilde{Z}(E[-k]))$ and some point $y$ in $\bP^{1}$,
   $E'$ is isomorphic to the object $\cO_{y}$.  Here, for some object
   $E'$ in $\tilde{Z}^{-1}(\tilde{Z}(E[-k]))$, we suppose otherwise and
   show that for any point $x$ in $\bP^{1}$ and any integer $i$, \
   $\Hom_{\cT}^{i}(E', \cO_{x})$ is the zero vector space.  For any
   point $x$ in $\bP^{1}$, since the objects $E'$ and $\cO_{x}$ are
   stable with the same phases, $\Hom_{\cT}(E', \cO_{x})$ and
   $\Hom_{\cT}(\cO_{x}, E')$ are the zero vector spaces; so,
   $\Hom_{\cT}^{2}(E', \cO_{x}) \cong \Hom_{\cT}(\cO_{x}, E')^{*}$ is
   also the zero vector space.  Since the objects $E'$ and $\cO_{x}$ are
   of the heart of $\sigma$, for any negative integer $i$, \
   $\Hom_{\cT}^{i}(E', \cO_{x})$ is the zero vector space, and also for
   any integer $i>2$, \ $\Hom_{\cT}^{i}(E', \cO_{x}) \cong
   \Hom_{\cT}^{2-i}(\cO_{x}, E')^{*}$ is the zero vector space. Now,
   since $\sigma$ is faithful, $\chi(\cO_{x}, \cO_{x}) =\chi(E',
   \cO_{x})=-\dim \Hom^{1}_{\cT}(E', \cO_{x})$ is zero.

   For the other case, we apply the twist functor
   $T_{\cO_{\bP^{1}}(w-1)}$ and use the above argument; by \cite[Lemma
   4.15 (i)(1)]{IshUeh}, \
   $T_{\cO_{\bP^{1}}(w-1)}(\cO_{\bP^{1}}(w-1)[1])=\cO_{\bP^{1}}(w-1)$
   and $T_{\cO_{\bP^{1}}(w-1)}(\cO_{\bP^{1}}(w))=\cO_{\bP^{1}}(w-2)[1]$.
  \end{proof}
  
  Before discussing the other cases, let us introduce some notions.  For
  each semistable object $E$ and a Jordan-H\"older decomposition
  $E\supset E_{1}\supset \cdots \supset E_{n}\supset 0$ in
  $\cP(\phi(E))$, {\it the cycle of simple components} \cite[Section
  2]{Ses} associated to the decomposition is defined to be $E_{n}\oplus
  E_{n-1}/E_{n-2}\oplus \cdots \oplus E/E_{1}$, which is by
  \cite[Theorem 2.1]{Ses}, up to isomorphisms, independent of the
  choices of Jordan-H\"older decompositions.  Semistable objects with
  isomorphic cycles of simple components are said to be {\it
  S-equivalent} \cite[Remark 2.1]{Ses}, \cite[Section 0]{MR466475}.
   
   Now, for some integer $|n|>1$, a semistable object $E$ with
   $[E]=n[\cO_{x}]$ in $K(\cT)$, the integer $k=\phi(E)-\phi(\cO_{x})$,
   and any object $E'$ of $\tilde{Z}^{-1}(\tilde{Z}(E[-k]))$, since for
   some point $x$ in $\bP^{1}$, \ $\Hom_{\cT}(E',\cO_{x})$ is not zero,
   $E'$ has a Jordan-H\"older decomposition whose composition factors
   are $\cO_{x}$ for points $x$ in $\bP^{1}$. So, up to S-equivalence on
   $\tilde{Z}^{-1}(\tilde{Z}(E))$ and autoequivalences on $\Stab(\cT)$,
   for each object of $\tilde{Z}^{-1}(\tilde{Z}(E))$, we have the
   corresponding $n$-fold direct sum of points in $\bP^{1}$.

\end{document}